\documentclass[10pt]{article}
\usepackage{amsfonts}
\usepackage{amssymb}
\usepackage{amsmath}
\usepackage{amsthm}
\usepackage{latexsym}
\usepackage{graphicx}
\usepackage{booktabs}
\usepackage{needspace}
\usepackage{color}
\newcommand{\V}[1]{\boldsymbol{#1}}

\newcommand{\Section}[1]{
   \Needspace{7\baselineskip}
   \par
   \refstepcounter{section}
   \bigskip\noindent
   {\large\bfseries\hyphenpenalty=10000 \exhyphenpenalty=10000
    \hbox{\thesection~~}#1\par}
   \nopagebreak
   \medskip
   \renewcommand{\theequation}{\thesection.\arabic{equation}}
   \setcounter{equation}{0}
   \setcounter{subsection}{0}
}

\newtheorem{thm}{Theorem}[section]
\newtheorem{lem}{Lemma}[section]
\newtheorem{rem}{Remark}[section]
\newtheorem{pro}{Proposition}[section]

\newtheorem{asm}{Assumption}[section]
\newtheorem{cor}{Corollary}[section]

\newcommand{\reff}[1]{(\ref{#1})}
\newcommand{\Proof}{\noindent {\bf Proof.~}\ }
\newcommand{\Sketch}{\noindent {\bf Sketch of proof.~}\ }
\newcommand{\Endproof}{$\hfill\Box$}
\newcommand{\Thm}[1]{Theorem #1}
\newcommand{\Lem}[1]{Lemma #1}
\newcommand{\Rem}[1]{Remark #1}
\newcommand{\Pro}[1]{Proposition #1}
\newcommand{\Cor}[1]{Corollary #1}
\newcommand{\Asm}[1]{Assumption #1}
\newcommand{\Tab}[1]{Table #1}
\newcommand{\Fig}[1]{Figure #1}

\newcommand{\be}{\begin{equation}}
\newcommand{\ee}{\end{equation}}
\newcommand{\p}{\partial}

\newcommand{\norm}[1]{\left\Vert#1\right\Vert}
\newcommand{\norme}[1]{\left\Vert{\hskip -2.6pt}\left\vert #1
     \right\vert{\hskip -2.6pt}\right\Vert}
\newcommand{\abs}[1]{\left\vert#1\right\vert}
\newcommand{\pd}[1]{\left\langle #1\right\rangle}
\newcommand{\set}[1]{\left\{#1\right\}}
\newcommand{\av}[1]{\left\{#1\right\}}
\newcommand{\jm}[1]{\left[#1\right]}

\newcommand{\mR}{\mathbb{R}}
\newcommand{\mP}{\mathbb{P}}
\newcommand{\Th}{{\cal T}_h}
\newcommand{\ThG}{{\cal T}_h^\Gamma}
\newcommand{\Fh}{{\cal F}_h}
\newcommand{\Gh}{\Gamma_h}
\newcommand{\Oh}[1]{\Omega_{#1}^h}
\newcommand{\Omh}{\Omega_h}
\newcommand{\nh}{\V{n}_h}
\newcommand{\nst}{\V{n}^*}
\newcommand{\rh}{\varrho_h}
\newcommand{\ph}{p_h}
\newcommand{\delh}{\delta_h}
\newcommand{\thh}{\theta_h}
\newcommand{\Tm}{{\cal T}_m}
\newcommand{\Lm}{\Lambda_m}
\newcommand{\Pim}{\Pi_m}
\newcommand{\jh}{j_h}
\newcommand{\Vh}{V_h}
\newcommand{\kbar}{\overline{\kappa}}
\newcommand{\thg}{\theta_\Gamma}
\newcommand{\sg}{s_\Gamma}

\begin{document}
\centerline{\Large \bf A Cut Finite Element Method with Boundary Value}

\vskip .2cm
\centerline{{\Large \bf Correction for Robin, Neumann and}}

\vskip .2cm
\centerline{{\Large \bf Interface Conditions}}

\vskip .5cm
\centerline{{\bf Erik Burman}$^{1}$\quad{\bf Peter Hansbo}$^{2}$\quad{\bf
Mats G. Larson}$^{3}$}

\vskip .3cm
\centerline{\it $^{1}$Department of Mathematics, University College London,}
\centerline{\it Gower Street, London WC1E 6BT, UK.}

\vskip .2cm
\centerline{\it $^{2}$Department of Mechanical Engineering, J\"onk\"oping
University,}
\centerline{\it SE-551 11 J\"onk\"oping, Sweden.}

\vskip .2cm
\centerline{\it $^{3}$Department of Mathematics and Mathematical Statistics,}
\centerline{\it Ume{\aa} University, SE-901 87 Ume{\aa}, Sweden.}

\vskip .3cm
\centerline{\small E-mail addresses: e.burman@ucl.ac.uk,
peter.hansbo@ju.se, mats.larson@umu.se.}

\vskip .5cm \noindent\rule[2mm]{\textwidth}{.1pt}

\noindent{\bf Abstract.}  We extend the cut finite element method with
boundary value correction of Burman, Hansbo and Larson [Math.\ Comp.\ 87
(2018) 633--657] to Robin and Neumann conditions and to interface problems
with imperfect bonding.  The formulation retains a piecewise affine
boundary for all integration, while using the exact geometry only at
boundary quadrature points.

The correction uses Taylor expansions of the solution and its gradient
along the discrete normal, with the exact normal used to evaluate the
transported flux.  At an interface, a further
correction to the normal flux jump accounts for the difference between
these two normals.  Under the stated trace and Poincar\'e bounds, the
corrected form is an $O(h^{1/2})$ perturbation of a symmetric form, whose
matrix can therefore be used as a preconditioner.

For a fixed finite Robin coefficient, we obtain the energy error bound
$C_\beta(h^p+\delh^{m+1})$, whereas the bound uniform in the Robin
coefficient is $C(h^p+h^{-1/2}\delh^{m+1})$, provided the Poincar\'e
constant is uniform.  With $\delh=O(h^2)$, the sufficient Taylor orders
are therefore $2m+2\ge p$ and $2m+\frac32\ge p$, respectively.
Numerical examples for $p=1,2,3$ illustrate these rates and the effect of
the flux-jump correction in the interface problem.


\noindent\rule[2mm]{\textwidth}{.1pt}

\Section{Introduction.}
\label{sec:intro}

Unfitted finite element methods use meshes that need not align with the
boundary or interface of the physical domain, which considerably
simplifies mesh generation.  The geometric approximation must nevertheless
be sufficiently accurate to support the desired convergence rate.  If the
boundary is replaced by a piecewise affine surface $\Gh$, imposing the
boundary condition directly on $\Gh$ limits this rate, regardless of the
polynomial degree.

There are two main approaches to this problem.  A higher order geometric
representation can be used, with curved elements or quadrature on the
exact domain.  Alternatively, the affine approximation can be retained and
the boundary values corrected through a Taylor expansion from $\Gamma$
to $\Gh$.  This idea was introduced by Bramble, Dupont and Thom\'ee
\cite{BDT72} for polygonal approximations of curved domains.  Burman,
Hansbo and Larson \cite{BHL18} subsequently developed it for unfitted
methods with Nitsche boundary conditions and ghost penalty stabilisation.
Their analysis shows that a correction of order $m$ gives the optimal
energy error for elements of degree $p$ when $2m+\frac32\ge p$.
Thus one Taylor term suffices for $p=2,3$ and two for
$p=4,5$ \cite[Tab.~1]{BHL18}.  Without correction, the geometric distance
must satisfy $\delh\lesssim h^{p+1/2}$ for the energy estimate and
$\delh\lesssim h^{p+1}$ for the $L^2$ estimate
\cite[Rem.~3.5]{BHL18}.  With one correction term, a piecewise affine
boundary with $\delh=O(h^2)$ suffices for degrees up to three.

The shifted boundary method of Main and Scovazzi \cite{MS18a,MS18b} uses
the same principle, with a surrogate boundary formed by mesh faces.
The computation then takes place on the resulting fitted submesh, and we
refer to \cite{ACS20} for the analysis.  Related methods on polytopal
approximations were developed in \cite{CPBG19,DGS25}, while the virtual
element methods in \cite{BPP22,BMPP24} also allow staircase boundaries
from pixel tessellations.  For unfitted meshes, boundary value correction
has also been combined with Lagrange multipliers \cite{BHL20}, a
penalty-free Nitsche method \cite{BBCL18}, and a mixed formulation
\cite{HLZ26}.  An alternative based on isoparametric mappings of cut
elements is given in \cite{L16}.  These developments mainly concern
Dirichlet boundary conditions.

For Neumann conditions, Bertoluzza, Montardini and Prada
\cite{bertoluzza2026highordervirtualelement} develop a high order virtual
element method on approximate domains.  Their formulation allows
non-aligned physical and numerical normals, with optimal error estimates
under suitable assumptions for arbitrary approximation order.  This
provides a related treatment of the geometric difficulties that arise
when the boundary condition involves the flux.

A Dirichlet condition prescribes the solution trace, so its transport
requires only a Taylor expansion of that trace.  The normal flux enters
through integration by parts on $\Omh$, where the resulting terms use the
normal of $\Gh$ and are exact on the discrete domain.

For a Robin condition, the boundary data also involve the flux, and we
consider
\be
  \label{eq:robin-intro}
  \kappa\,\p_{\V n}u+\beta u=g\qquad\hbox{on }\Gamma ,
\ee
where $\V n$ is the exact normal.  Since the discrete normal $\nh$
generally differs from $\V n$ by $O(h)$, transporting the evaluation point
while retaining $\nh$ leaves an $O(h)$ error in the flux.  We remove
this error by transporting the vector field $\kappa\nabla u$ and then
taking its scalar product with the exact normal $\nst$ at the image
point.  The resulting operator, denoted by $\Lm$, has an error of
$O(\delh^{m+1})$, as for the trace correction
(\Lem{\ref{lem:rem}}).  The construction uses the polynomial gradient on
each cut element together with pointwise evaluations of the exact normal,
without requiring curved quadrature or additional unknowns.

We combine this flux correction with a formulation that is uniform in the
Robin coefficient $\beta$, connecting the Neumann problem at $\beta=0$
to the Dirichlet limit as $\beta\to\infty$.  Following
Juntunen and Stenberg \cite{JS09}, we set $\rho:=\gamma\kappa/h$ and
introduce the weights
\be
  \label{eq:theta-intro}
  \theta:=\frac{\beta}{\beta+\rho}\in[0,1],\qquad
  s:=\frac{1}{\beta+\rho} ,
\ee
which give a single formulation for the whole family.  At $\theta=1$ and
$s=0$, it reduces to the Dirichlet method with boundary value correction
from \cite{BHL18}, while $\theta=0$ gives a corrected Neumann formulation.
The pair $(\theta,s)$ is equivalent to the weights in
\cite[(2.5)--(2.6)]{JS09}, but makes the roles of the trace and flux
corrections explicit.  The estimates have constants independent
of $\beta$ under the stated coercivity and regularity assumptions.

We also consider an elliptic interface problem, where the prescribed flux
jump is taken with respect to the exact normal,
$\jm{\kappa\p_{\V n}u}=g_N$.  Integration by parts on the two discrete
subdomains instead gives the jump with respect to $\nh$.  Their
difference contains
\be
  \label{eq:mis-intro}
  \jm{\kappa\nabla u}\cdot(\nh-\nst).
\ee
This term depends on the tangential flux jump, which is not prescribed by
the interface data and therefore cannot be removed by transporting those
data alone.  We retain the corresponding integration-by-parts term,
$-\pd{\jm{\kappa\p_{\nh}u},\av{v}'}_{\Gh}$, in the formulation on
$\Gh$.  Although this term vanishes on the exact interface for a
homogeneous normal flux jump, it need not vanish on the surrogate
interface.  We control it using two different trace estimates
(\Lem{\ref{lem:mis}}).  The elementwise gradient requires a discrete
trace estimate with a factor $h^{-1/2}$, whereas the test function belongs
to $H^1$ and satisfies a continuous trace estimate without this factor.
Together with $\norm{\nh-\nst}=O(h)$, this gives an $O(h^{1/2})$
perturbation.

Finally, we consider imperfect bonding, where a spring layer with
stiffness $\beta_\Gamma$ joins the two phases, as in \cite{HH04,HH02}.
For this interface counterpart of the Robin condition, the weight
$\thg=\beta_\Gamma/(\beta_\Gamma+\rho)$ gives a single formulation for
all bonding stiffnesses.  Uniform estimates additionally require the
phasewise Poincar\'e bounds stated below.  In \S\ref{sec:HH}, we interpret
the boundary value correction as an artificial compliance of order $h^2$,
given by the geometric offset weighted by the jump of the inverse
diffusivity.  This correction adds to the compliance in the physical
bonding law.

We give the main arguments in \Lem{\ref{lem:rem}},
\Lem{\ref{lem:mis}}, \Thm{\ref{thm:cons}} and
\Thm{\ref{thm:icons}}, while proofs that follow from standard arguments
are sketched.

Throughout, $a\lesssim b$ means $a\le Cb$, where $C$ is independent of
the mesh size, the position of the boundary relative to the mesh,
$\beta$, and the contrast, unless stated otherwise.  It may depend on the
polynomial degree, the shape regularity, the curvature of $\Gamma$, and
the stated Poincar\'e bounds.  Uniformity in the physical parameters is
therefore conditional on uniform control of these bounds and of the
solution norms.

\Section{The Robin problem and its two limits.}
\label{sec:model}

Let $\Omega\subset\mR^d$, $d=2,3$, be a bounded domain with a $C^2$
boundary $\Gamma=\p\Omega$ and outward unit normal $\V n$.  For constant
$\kappa>0$, reaction coefficient $c\ge0$, and given $f$, $g$, and
$\beta\in[0,\infty]$, we consider
\be
  \label{eq:model}
  \begin{array}{rcll}
   -\nabla\cdot(\kappa\nabla u)+cu&=&f&\hbox{in }\Omega,\\[2pt]
   \kappa\,\p_{\V n}u+\beta u&=&g&\hbox{on }\Gamma .
  \end{array}
\ee
The case $\beta=0$ gives the Neumann problem, while finite positive
$\beta$ gives the Robin problem.  The Dirichlet limit is obtained by setting
$g=\beta g_D$ and letting $\beta\to\infty$ with $g_D$ fixed.  We assume
$c\ge c_0>0$ or $\beta\ge\beta_0>0$.
For the pure Neumann problem ($c=\beta=0$), the continuous solution is
defined modulo constants and the data must satisfy the usual compatibility
condition.  We restrict both the discrete trial and test spaces to functions
with zero mean on $\Omh$, and normalise the exact solution to have the
same mean on $\Omh$ in the error analysis.  This restriction is needed
because flux transport need not preserve the discrete compatibility
condition for constant test functions.
In \cite{JS09}, the condition is written in terms of the compliance
$\epsilon=1/\beta$.  We use $\beta$ so that the Neumann condition
corresponds to the finite value $\beta=0$ and the notation agrees with the
bonding stiffness in \S\ref{sec:iface}.

We denote the normal flux by
\be
  \label{eq:lambda}
  \lambda(u):=\kappa\,\p_{\V n}u=\kappa\,\V n\cdot\nabla u
  \qquad\hbox{on }\Gamma .
\ee
The boundary condition then reads $\beta u+\lambda(u)=g$, while
multiplication of the equation by a test function and integration by parts
give
\be
  \label{eq:ibp}
  (\kappa\nabla u,\nabla v)_\Omega+(cu,v)_\Omega
     -\pd{\lambda(u),v}_\Gamma=(f,v)_\Omega .
\ee
The method combines this identity with the boundary condition through the
following substitution.

\begin{lem}[The uniform substitution]
\label{lem:subst}
Let $\rho>0$, let $\beta\in[0,\infty)$, and let $\theta$ and $s$ be given
by \reff{eq:theta-intro}.  Then $\theta=\beta s$ and
$1-\theta=\rho s$.  For any $U$, $\Lambda$, and $G$ satisfying
$\beta U+\Lambda=G$, we have
\be
  \label{eq:subst}
  \Lambda=\rho s\,G-\theta\rho\,U+\theta\Lambda,
  \qquad
  \theta U+s\Lambda=s\,G .
\ee
\end{lem}

\Proof
The first identity follows from
$(1-\theta)\Lambda=\rho s\Lambda=\rho s(G-\beta U)=\rho sG-\theta\rho U$.
The second is $s(\beta U+\Lambda)=sG$. \Endproof

The first identity replaces part of the flux in \reff{eq:ibp} by the
trace and the data.  The second supplies a residual that may be tested
without changing the exact solution.  At $\beta=0$, the identities become
$\Lambda=G$ and $s\Lambda=sG$, while in the Dirichlet limit they become
$\Lambda=\Lambda$ and $U=g_D$.

\begin{rem}[Relation to the method of Juntunen and Stenberg]
\label{rem:JS}
Applying \Lem{\ref{lem:subst}} on the exact domain with $U=u$ and
$\Lambda=\lambda(u)$, without geometric correction, gives the symmetric
formulation
$(\kappa\nabla u,\nabla v)_\Omega+(cu,v)_\Omega
 -\theta(\pd{\lambda(u),v}_\Gamma+\pd{u,\lambda(v)}_\Gamma)
 +\theta\rho\pd{u,v}_\Gamma-s\pd{\lambda(u),\lambda(v)}_\Gamma
 =(f,v)_\Omega+\rho s\pd{g,v}_\Gamma-s\pd{g,\lambda(v)}_\Gamma$.
It agrees with \cite[(2.5)--(2.6)]{JS09} under the correspondence
$\epsilon=1/\beta$ and $\gamma_{\rm JS}h_E=1/\rho$.
Their stability condition $\gamma_{\rm JS}<1/C_I$
\cite[Thm.~3.2]{JS09} requires sufficiently large $\gamma$ in our
notation.  We extend this form by introducing a surrogate boundary, a
ghost penalty, and the two transport operators.
\end{rem}

\Section{The discrete geometry and the two Taylor operators.}
\label{sec:geom}

\subsection{The surrogate boundary.}

Let $\Th$ be a quasi-uniform, shape-regular simplicial mesh of a fixed
polytopal domain containing $\overline\Omega$, with mesh size
$h=\max_{T\in\Th}h_T$.  We describe the boundary by a level set function
$\Phi\in C^2$, with $\abs{\nabla\Phi}>0$ near $\Gamma$, such that
\be
  \label{eq:ls}
  \Omega=\set{\Phi<0},\qquad\Gamma=\set{\Phi=0},\qquad
  \V n=\frac{\nabla\Phi}{\abs{\nabla\Phi}}\quad\hbox{on }\Gamma .
\ee
The cut elements are
$\ThG:=\set{T\in\Th:\ T\cap\Gamma\ne\emptyset}$.  In two dimensions,
we assume that each cut element meets $\Gamma$ in one connected arc with
endpoints on two element edges, and replace this arc by the chord
$\Gh^T$ joining its endpoints.  In three dimensions we use a piecewise
planar surface and assume the geometric estimates stated below.  We set
\be
  \label{eq:Gh}
  \Gh:=\bigcup_{T\in\ThG}\Gh^T .
\ee
We assume that $\Gh$ is a closed Lipschitz surface bounding a polytopal
domain $\Omh$, with uniformly bounded Lipschitz character as $h\to0$.
The active mesh $\Th^{\rm act}$ consists of the elements for which
$\abs{T\cap\Omh}>0$, and on their union $\Omega_h^{\rm act}$ we define
\be
  \label{eq:space}
  \Vh:=\set{v\in C^0(\Omega_h^{\rm act}):\ v|_T\in\mP_p(T)\ \
     \forall T\in\Th^{\rm act}} .
\ee
The outward unit normal of $\Omh$, denoted by $\nh$, is constant on each
affine piece of $\Gh$.  Since all integrals are taken over $\Omh$, $\Gh$,
or faces of the active mesh, the method requires quadrature only on
polytopal sets.

\subsection{The three pointwise maps.}

Let $x\in\Gh$ be a point where $\nh$ is defined.  For sufficiently small
$h$, the equation $\Phi(x+r\,\nh(x))=0$ has a unique root near zero.  We
define
\be
  \label{eq:maps}
  \rh(x):=\hbox{that root},\qquad
  \ph(x):=x+\rh(x)\,\nh(x)\in\Gamma,\qquad
  \nst(x):=\frac{\nabla\Phi(\ph(x))}{\abs{\nabla\Phi(\ph(x))}} .
\ee
Thus $\ph$ transfers points from $\Gh$ to $\Gamma$ along the discrete
normal, while $\nst(x)$ gives the exact normal at the image point
$\ph(x)$.  The geometric errors are measured by
\be
  \label{eq:dhth}
  \delh:=\norm{\rh}_{L^\infty(\Gh)},\qquad
  \thh:=\norm{\nh-\nst}_{L^\infty(\Gh)} .
\ee

\begin{lem}[Geometric consistency]
\label{lem:geom}
In two dimensions, let $\Gamma$ be of class $C^2$ with curvature bounded
by $c_\Gamma$.  For sufficiently small $c_\Gamma h$, the chord
construction satisfies
\be
  \label{eq:geomest}
  \delh\le c_\Gamma h^2,\qquad
  \thh\le 2c_\Gamma h,\qquad
  \abs{\Omh\triangle\Omega}\le c_\Gamma\abs{\Gamma}h^2 .
\ee
The surface Jacobian $J_{\ph}$ of $\ph$ satisfies
$\abs{J_{\ph}-1}\lesssim h^2$ on each affine piece of $\Gh$.
In three dimensions, these estimates are assumptions on the piecewise
planar approximation.
\end{lem}

\Proof
For a chord of length $\ell\le h_T$ with endpoints on $\Gamma$, consider
a second-order expansion of the corresponding arc.  This gives a distance
of order $c_\Gamma\ell^2$ and a normal error of order $c_\Gamma\ell$,
from which the first two estimates follow for sufficiently small $h$.
Integration of the distance estimate over $\Gh$ then gives
the bound on the symmetric difference.  On each chord, the surface
Jacobian differs from one by
$O(c_\Gamma\delh+\thh^2)$, which is $O(h^2)$ with a constant depending
on the curvature. \Endproof

Other expansion directions are possible, and \cite{BHL18} allows a
field $\V\nu_h$ satisfying $\norm{\V\nu_h-\V n\circ p}=O(1)$.
With our choice $\V\nu_h=\nh$, the offset $\rh$ can be computed by
Newton iteration for
$\Phi(x+r\nh(x))=0$, starting at $r=0$.  Since $\ph$ is not the closest
point map, we evaluate $\nst$ separately from $\nabla\Phi$ for use in
the flux transport \reff{eq:Lm}.  All three maps are needed only at
quadrature nodes on $\Gh$.

\subsection{Transporting a trace and transporting a flux.}

Write $\p_{\nh}:=\nh\cdot\nabla$ for differentiation in the discrete
normal direction.  For $m\ge0$ and a function $w$ smooth near $x\in\Gh$,
define
\be
  \label{eq:Tm}
  \Tm(w)(x):=\sum_{j=0}^{m}\frac{\rh(x)^j}{j!}\,\p_{\nh}^{\,j}w(x) .
\ee
This is the Taylor polynomial along the segment from $x$ to $\ph(x)$,
evaluated at $\ph(x)$, with ${\cal T}_0=\hbox{Id}$.
The operator transports the trace, as in \cite{BDT72,BHL18}, and we also
apply it componentwise to vector fields.

To transport the flux in the Robin condition, we define
\be
  \label{eq:Lm}
  \Lm(w)(x):=\Tm\big(\kappa\nabla w\big)(x)\cdot\nst(x)
   =\kappa\sum_{j=0}^{m}\frac{\rh(x)^j}{j!}\,
     \big(\p_{\nh}^{\,j}\nabla w\big)(x)\cdot\nst(x) .
\ee
We first expand $\kappa\nabla w$ along $\nh$ and then take its scalar
product with the exact normal at the image point.  Since $\nst(x)$ is a
fixed vector at each quadrature node, this is equivalent to expanding the
scalar function $y\mapsto\kappa\,\nst(x)\cdot\nabla w(y)$.
In particular, no derivative of the normal is required.  For a discrete
function, $\Tm$ uses derivatives up to order $m$ and $\Lm$ uses
derivatives up to order $m+1$, evaluated from the polynomial on the cut
element.  We therefore take $p\ge m+1$.

\begin{lem}[Transport error]
\label{lem:rem}
Let $w\in W^{m+2,\infty}(U)$, where $U$ is a neighbourhood of the segment
$[x,\ph(x)]$.  Then
\be
  \label{eq:rem}
  \begin{array}{l}
  \displaystyle
  \abs{\Tm(w)(x)-w(\ph(x))}\le
   \frac{\delh^{m+1}}{(m+1)!}\,\norm{\p_{\nh}^{m+1}w}_{L^\infty},
  \\[10pt]\displaystyle
  \abs{\Lm(w)(x)-\lambda(w)(\ph(x))}\le
   \frac{\kappa\,\delh^{m+1}}{(m+1)!}\,\norm{\p_{\nh}^{m+1}\nabla w}_{L^\infty} .
  \end{array}
\ee
\end{lem}

\Proof
The first estimate follows by applying Taylor's theorem with remainder to
$t\mapsto w(x+t\nh(x))$ between $0$ and $\rh(x)$.
For the second, apply the same argument to
$t\mapsto\kappa\,\nst(x)\cdot\nabla w(x+t\nh(x))$.
Its value at $t=\rh(x)$ is $\lambda(w)(\ph(x))$, since $\nst(x)$ is the
exact normal at $\ph(x)$. \Endproof

The discrete flux $\kappa\p_{\nh}w(x)$ differs from the exact flux at the
image point by $\kappa\nabla w\cdot(\nh-\nst)+O(\delh)$, which is
generally $O(h)$ for an affine boundary approximation.  A Taylor expansion
that retains $\nh$ in the scalar product leaves this directional error.
By using $\nst$ in $\Lm$, we retain only the transport error
$O(\delh^{m+1})$.  This construction uses the elementwise polynomial
gradient and does not require gradient continuity between elements.

\begin{asm}[Extension]
\label{asm:ext}
The solution $u$ of \reff{eq:model} and the datum $f$ extend to a fixed
neighbourhood of $\overline{\Omega}$, with
$u\in W^{m+2,\infty}\cap H^{p+1}$ there.  The extensions satisfy
$-\nabla\cdot(\kappa\nabla u)+cu=f$ pointwise in
$\Omh\cup\Omega$.  The datum $g$ also extends to a neighbourhood of
$\Gamma$.
\end{asm}

This assumption allows us to integrate by parts on $\Omh$ without a
volume consistency term.  Under the weaker extension assumption in
\cite{BHL18}, the residual outside $\Omega$ is retained and estimated on
the strip $\Omh\setminus\Omega$ \cite[(3.14), (3.49)]{BHL18}.
We use \Asm{\ref{asm:ext}}, which holds exactly for all test cases in
\S\ref{sec:num}, to focus on the transport of the boundary condition.

\Section{The method.}
\label{sec:method}

\subsection{The form.}

Let $\gamma>0$ be a Nitsche parameter and set
\be
  \label{eq:rho}
  \rho:=\frac{\gamma\kappa}{h},\qquad
  \theta:=\frac{\beta}{\beta+\rho},\qquad s:=\frac{1}{\beta+\rho} .
\ee
We use a ghost penalty on the interior faces of the active mesh adjacent
to cut elements,
\be
  \label{eq:ghost}
  \jh(v,w):=\sum_{F\in\Fh^\Gamma}\sum_{j=1}^{p}
     \gamma_g\,\kappa\,h^{2j-1}\pd{\jm{\p_{\V n_F}^{\,j}v},
                                   \jm{\p_{\V n_F}^{\,j}w}}_F .
\ee
The discrete problem is to find $u_h\in\Vh$ such that
$A_h(u_h,v)=L_h(v)$ for all $v\in\Vh$, where
\be
  \label{eq:Ah}
  \begin{array}{rcl}
  A_h(u,v)&:=&(\kappa\nabla u,\nabla v)_{\Omh}+(cu,v)_{\Omh}+\jh(u,v)
      -\pd{\kappa\p_{\nh}u,\,v}_{\Gh}\\[4pt]
   &&{}+\ \rho s\pd{\Lm(u),v}_{\Gh}
       +\theta\rho\pd{\Tm(u),v}_{\Gh}\\[4pt]
   &&{}-\ \theta\pd{\Tm(u),\kappa\p_{\nh}v}_{\Gh}
       -s\pd{\Lm(u),\kappa\p_{\nh}v}_{\Gh},
  \end{array}
\ee
\be
  \label{eq:Lh}
  L_h(v):=(f,v)_{\Omh}+\rho s\pd{g\circ\ph,\,v}_{\Gh}
                      -s\pd{g\circ\ph,\,\kappa\p_{\nh}v}_{\Gh} .
\ee
If $c=\beta=0$, both trial and test functions are restricted to the
subspace $\set{v\in\Vh:(v,1)_{\Omh}=0}$.
The datum is evaluated at $\ph(x)\in\Gamma$, while $\Tm$ and $\Lm$
approximate the corresponding solution trace and flux from the polynomial
on the cut element.

The term $-\pd{\kappa\p_{\nh}u,v}_{\Gh}$ comes from integration by parts
on $\Omh$ and therefore uses the normal $\nh$ of that domain.
It is present for every $\beta$, while the terms involving $\Tm$ and
$\Lm$ impose the transported boundary condition through
\Lem{\ref{lem:subst}}.  The ghost penalty provides uniform control of the
active mesh norm with respect to the boundary cuts.

\subsection{The two limits and the Taylor order.}

Let $g=\beta g_D$ and keep $g_D$ fixed as $\beta\to\infty$.
Then $\theta\to1$, $s\to0$, $\rho s\,g\to\rho g_D$, and
$s\,g\to g_D$.  The method tends to Nitsche's method with the boundary
value correction of \cite{BHL18}, acting on the trial function only.
The flux transport $\Lm$ disappears because the Dirichlet condition
prescribes only the trace.

At $\beta=0$, we have $\theta=0$, $s=1/\rho$, and $\rho s=1$.
Equations \reff{eq:Ah}--\reff{eq:Lh} become
\[
\begin{aligned}
 &(\kappa\nabla u,\nabla v)_{\Omh}+(cu,v)_{\Omh}+\jh(u,v)
 -\pd{\kappa\p_{\nh}u-\Lm(u),v}_{\Gh}\\
 &\qquad-\frac1\rho\pd{\Lm(u),\kappa\p_{\nh}v}_{\Gh}
 =(f,v)_{\Omh}+\pd{g\circ\ph,v}_{\Gh}
 -\frac1\rho\pd{g\circ\ph,\kappa\p_{\nh}v}_{\Gh}.
\end{aligned}
\]
The two terms with coefficient $1/\rho$ form a Barbosa--Hughes
stabilisation whose residual vanishes at the exact solution on the exact
geometry, as in \cite[(2.14)]{JS09}.  On $\Gh$, the residual instead has
the Taylor transport error from \Lem{\ref{lem:rem}}.  Dropping these two terms gives
a natural Neumann formulation with a flux correction: the flux from
integration by parts is replaced by $\Lm(u)$ before the datum is imposed.

\begin{rem}[Where the Robin coefficient acts]
\label{rem:where}
The trace correction $\Tm$ is multiplied by $\theta$, while the flux
correction $\Lm$ is multiplied by $\rho s=1-\theta$ or by
$s=(1-\theta)/\rho$.  Thus the Dirichlet limit uses only the trace
correction and the Neumann case only the flux correction, with both
entering for finite positive $\beta$.
\end{rem}

The uniform consistency bound in \Cor{\ref{cor:cons}} is
$O(h^{-1/2}\delh^{m+1})$, which becomes $O(h^{2m+3/2})$ under
$\delh\lesssim c_\Gamma h^2$.  It is therefore bounded by the
approximation error $O(h^p)$ when
\be
  \label{eq:rule}
  2m+\tfrac32\ \ge\ p ,
  \qquad\hbox{i.e.}\qquad
  m=0\ (p=1),\quad m=1\ (p=2,3),\quad m=2\ (p=4,5).
\ee
This is the rule in \cite[Tab.~1]{BHL18} for the Dirichlet problem, and
also suffices for the expected $L^2$ estimate
(\Rem{\ref{pro:l2}} and \cite[Tab.~2]{BHL18}).
The estimate can be sharpened when $\beta$ is fixed and finite, as shown
below.

\Section{Analysis.}
\label{sec:anal}

\subsection{Norms and trace inequalities.}

We use the energy norm
\be
  \label{eq:norm}
  \norme{v}_h^2:=\norm{\sqrt\kappa\nabla v}^2_{\Omh}+c\norm{v}^2_{\Omh}
    +\jh(v,v)+\theta\rho\norm{v}^2_{\Gh} .
\ee
The boundary weight is $\rho=\gamma\kappa/h$ at the Dirichlet limit and
vanishes at the Neumann limit.  For fixed $0<\beta<\infty$,
$\theta\rho=\beta\rho/(\beta+\rho)\to\beta$ as $h\to0$, so the norm
approaches the natural Robin energy norm.
When $c=\beta=0$, we restrict the discrete space to functions with zero
mean on $\Omh$.

The ghost penalty controls the gradient on the full active elements.
Combined with the polynomial inverse inequality, this gives
\be
  \label{eq:disc}
  \begin{split}
  \kappa^2\norm{\nabla v}^2_{\Gh}
  &\le \frac{C_I\kappa}{h}
     \Big(\norm{\sqrt\kappa\nabla v}^2_{\Omh}+\jh(v,v)\Big),\\
  \norm{v}^2_{\Gh}
  &\lesssim h^{-1}\norm{v}^2_{\Omh}
       +h\Big(\norm{\nabla v}^2_{\Omh}+\kappa^{-1}\jh(v,v)\Big),
  \qquad v\in\Vh .
  \end{split}
\ee
The constants are independent of the cut position, and analogous estimates
hold on active-mesh faces.  These are standard consequences of the ghost
penalty, for which we refer to \cite{B10,BHL18}.

For the trace of the function itself, we also use
\be
  \label{eq:ctrace}
  \norm{v}_{L^2(\Gh)}\le C_{\rm tr}\norm{v}_{H^1(\Omh)},
  \qquad v\in H^1(\Omh).
\ee
The constant $C_{\rm tr}$ is independent of $h$ because the domains
$\Omh$ have a uniform Lipschitz character.
For a resolved chordal approximation of a $C^2$ curve, neighbouring
segments meet at an angle $\pi+O(h)$, which allows local boundary charts
to be chosen with uniform bounds.  In three dimensions we assume the same
uniform domain property.

We make the Poincar\'e bound explicit:
\be
  \label{eq:poincare}
  \kappa\norm{v}^2_{\Omh}\le\Lambda^2\norme{v}_h^2,
  \qquad v\in\Vh .
\ee
We assume that $\Lambda$ is independent of $h$ and make its role in the
analysis explicit.  If $c\ge c_0>0$, we may take
$\Lambda^2=\kappa/c_0$.
For $\beta\ge\beta_0>0$, a Poincar\'e inequality with a boundary term
also gives a uniform bound for sufficiently small $h$.
For $c=\beta=0$, the zero-mean restriction gives the usual Poincar\'e
bound.  Since the constants below may depend on $\Lambda$, uniformity in
$\beta$ requires a uniform bound on $\Lambda$.
In particular, such a bound does not hold on the full space as
$\beta\downarrow0$ when $c=0$.

Combining \reff{eq:ctrace} and \reff{eq:poincare} gives
\be
  \label{eq:energytrace}
  \norm{v}_{\Gh}\le C_\Lambda\kappa^{-1/2}\norme{v}_h,
  \qquad C_\Lambda:=C_{\rm tr}(1+\Lambda^2)^{1/2}.
\ee
This estimate has no negative power of $h$ and will be used in both the
stability and consistency arguments.

\begin{lem}[Taylor operators on discrete functions]
\label{lem:tpol}
Let $0\le m\le p-1$ and assume $\delh/h$ is bounded.
Then
\be
  \label{eq:tpert}
  \begin{split}
  \norm{\Tm(v)-v}_{\Gh}
   &\lesssim\frac{\delh}{\sqrt{\kappa h}}\norme{v}_h,\\
  \norm{\Lm(v)-\kappa\nst\cdot\nabla v}_{\Gh}
   &\lesssim\frac{\delh}{h}\Big(\frac{\kappa}{h}\Big)^{1/2}
      \norme{v}_h .
  \end{split}
\ee
Both left-hand sides vanish when $m=0$.
\end{lem}

\Proof
On each full active element, we apply the polynomial inverse inequality
and use the ghost penalty to control the gradient.
The $j$th trace and flux terms are then bounded, respectively, by
\[
 C\delh^j h^{-(j-1)}(\kappa h)^{-1/2}\norme{v}_h
 \quad\hbox{and}\quad
 C(\delh/h)^j(\kappa/h)^{1/2}\norme{v}_h.
\]
Summing over $j=1,\ldots,m$ proves the result. \Endproof

\subsection{The normal mismatch.}

Integration by parts uses the normal $\nh$, whereas the boundary condition
uses $\nst$.  For $m=0$, the resulting flux difference is
\be
  \label{eq:E}
  E(v):=\kappa(\nst-\nh)\cdot\nabla v,
  \qquad
  \abs{E(v)}\le\thh\kappa\abs{\nabla v}.
\ee
The term $\pd{E(v),v}_{\Gh}$ occurs in the stability estimate.

\begin{lem}[Bound for the normal mismatch]
\label{lem:mis}
Under \reff{eq:poincare},
\be
  \label{eq:mis}
  \abs{\pd{E(v),v}_{\Gh}}
     \le Cc_\Gamma h^{1/2}\norme{v}_h^2,
  \qquad v\in\Vh .
\ee
The constant depends on $C_I$ and $C_\Lambda$ but is independent of $h$
and of the cut position.
\end{lem}

\Proof
By the Cauchy--Schwarz inequality,
\[
 \abs{\pd{E(v),v}_{\Gh}}
 \le\thh\kappa\norm{\nabla v}_{\Gh}\norm{v}_{\Gh}.
\]
We apply the discrete trace estimate \reff{eq:disc} to the gradient and
the continuous estimate \reff{eq:energytrace} to the function itself.
Together they give
\[
 \abs{\pd{E(v),v}_{\Gh}}
 \le C_I^{1/2}C_\Lambda\thh h^{-1/2}\norme{v}_h^2.
\]
The result follows from $\thh\le2c_\Gamma h$. \Endproof

\begin{rem}[The Neumann case]
\label{rem:nocirc}
The volume terms and \reff{eq:poincare} control the $H^1$ norm without a
boundary penalty.  Hence \reff{eq:energytrace} and the argument above also
apply when $\beta=0$.
If $c=0$, both trial and test functions must satisfy the zero-mean
restriction specified above.  The corrected form need not annihilate
constant test functions, so a quotient-space statement alone is not enough.
\end{rem}

\subsection{Stability.}

\begin{thm}[Stability]
\label{thm:stab}
Let $\gamma$ be sufficiently large and let $0\le m\le p-1$.
Under \reff{eq:poincare}, there are $h_0>0$ and $C>0$ such that
\be
  \label{eq:stab}
  A_h(v,v)\ge C\norme{v}_h^2,
  \qquad v\in\Vh,\quad h\le h_0.
\ee
The constants may depend on $\gamma$, the mesh regularity, the geometry
and $\Lambda$, but are independent of the cut position.
They are also independent of $\beta$ whenever $\Lambda$ is uniformly
bounded.
\end{thm}

\Proof
For $m=0$, we have $\Tm(v)=v$ and
$\Lm(v)=\kappa\p_{\nh}v+E(v)$.  Using $\rho s=1-\theta$, we can therefore
write the form as
\be
  \label{eq:stabsplit}
  \begin{split}
  A_h(v,v)={}&\norm{\sqrt\kappa\nabla v}^2_{\Omh}
       +c\norm{v}^2_{\Omh}+\jh(v,v)
       +\theta\rho\norm{v}^2_{\Gh}\\
   &-2\theta\pd{\kappa\p_{\nh}v,v}_{\Gh}
       -s\norm{\kappa\p_{\nh}v}^2_{\Gh}+{\cal R}(v),
  \end{split}
\ee
where
\[
 {\cal R}(v)=(1-\theta)\pd{E(v),v}_{\Gh}
       -s\pd{E(v),\kappa\p_{\nh}v}_{\Gh}.
\]
Young's inequality gives
\[
 2\theta\abs{\pd{\kappa\p_{\nh}v,v}_{\Gh}}
 \le\frac{2\theta}{\rho}\norm{\kappa\p_{\nh}v}^2_{\Gh}
     +\frac{\theta\rho}{2}\norm{v}^2_{\Gh}.
\]
Since $(2\theta/\rho)+s\le2/\rho$, the first estimate in
\reff{eq:disc} absorbs the flux terms when $\gamma>4C_I$.
The remaining positive terms control $\frac12\norme{v}_h^2$.

By \Lem{\ref{lem:mis}}, the first term in ${\cal R}$ is bounded by
$Cc_\Gamma h^{1/2}\norme{v}_h^2$.
The bound $s\le h/(\gamma\kappa)$ gives
$Cc_\Gamma h\norme{v}_h^2$ for the second term, so both can be absorbed
for small $h$.

For $m\ge1$, \Lem{\ref{lem:tpol}} bounds the trace-correction terms and
the flux-flux term by
$C(\delh/h)\norme{v}_h^2$.
Applying \reff{eq:energytrace} to the remaining flux-correction term gives
the bound $CC_\Lambda\delh h^{-3/2}\norme{v}_h^2$.
Since $\delh\lesssim h^2$, these terms can also be absorbed. \Endproof

\subsection{Consistency.}

Define the transported boundary residual by
\be
  \label{eq:R}
  R(w):=\beta\Tm(w)+\Lm(w)-g\circ\ph
  \qquad\hbox{on }\Gh .
\ee

\begin{thm}[Consistency]
\label{thm:cons}
For $0\le\beta<\infty$, let $u$ solve \reff{eq:model} and satisfy
\Asm{\ref{asm:ext}}.  Then
\be
  \label{eq:cons}
  A_h(u,v)-L_h(v)=\rho s\pd{R(u),v}_{\Gh}
       -s\pd{R(u),\kappa\p_{\nh}v}_{\Gh},
  \qquad v\in\Vh,
\ee
and
\be
  \label{eq:Rbound}
  \norm{R(u)}_{\Gh}\lesssim
   \frac{\delh^{m+1}}{(m+1)!}
   \Big(\beta\norm{\p^{m+1}_{\nh}u}_{L^\infty}
     +\kappa\norm{\p^{m+1}_{\nh}\nabla u}_{L^\infty}\Big)
     \abs{\Gh}^{1/2}.
\ee
\end{thm}

\Proof
The extension assumption allows integration by parts on $\Omh$, giving
\[
 (\kappa\nabla u,\nabla v)_{\Omh}+(cu,v)_{\Omh}
       -\pd{\kappa\p_{\nh}u,v}_{\Gh}=(f,v)_{\Omh}.
\]
Since the extension is smooth across mesh faces, $\jh(u,v)=0$.
Subtracting the load from \reff{eq:Ah} therefore leaves
\[
 \rho s\pd{\Lm(u)-g\circ\ph,v}_{\Gh}
 +\theta\rho\pd{\Tm(u),v}_{\Gh}
 -\pd{\theta\Tm(u)+s\Lm(u)-sg\circ\ph,\kappa\p_{\nh}v}_{\Gh}.
\]
Using $\theta=\beta s$ gives \reff{eq:cons}, while the boundary condition
at $\ph(x)$ gives
\[
 R(u)=\beta\big(\Tm(u)-u\circ\ph\big)
       +\Lm(u)-\lambda(u)\circ\ph.
\]
The residual bound follows by applying \Lem{\ref{lem:rem}} to these two
remainders. \Endproof

\begin{cor}[Weighted consistency]
\label{cor:cons}
Let $C_u$ denote the $W^{m+2,\infty}$ norm of the extension of $u$ and set
$C_\Lambda:=C_{\rm tr}(1+\Lambda^2)^{1/2}$.  Then
\be
  \label{eq:consbound}
  \abs{A_h(u,v)-L_h(v)}
  \lesssim \delh^{m+1} C_u\sqrt\kappa
  \left[\left(\frac\theta h\right)^{1/2}
       +(1-\theta)(C_\Lambda+h^{1/2})\right]\norme v_h .
\ee
The constant is independent of $h$, $\beta$, and the position of the cut,
with the dependence on the Poincar\'e estimate displayed through
$C_\Lambda$.
For each fixed finite $\beta$, the right-hand side is
$O(\delh^{m+1})\norme v_h$.
At the Dirichlet limit it is $O(h^{-1/2}\delh^{m+1})\norme v_h$.
\end{cor}

\Proof
Write
\[
 r_T:=\Tm(u)-u\circ\ph,
 \qquad r_\Lambda:=\Lm(u)-\lambda(u)\circ\ph .
\]
This gives the decomposition $R(u)=\beta r_T+r_\Lambda$, with the transport
estimates
$\norm{r_T}_{\Gh}\lesssim\delh^{m+1}C_u$ and
$\norm{r_\Lambda}_{\Gh}\lesssim\kappa\delh^{m+1}C_u$.
We estimate the two residuals separately, starting with the trace
residual.  The boundary term in the energy norm gives
\[
 \theta\rho\norm v_{\Gh}
 \le(\theta\rho)^{1/2}\norme v_h .
\]
This inequality also holds at $\theta=0$, when both sides vanish.
The discrete flux estimate gives
$\theta\norm{\kappa\p_{\nh}v}_{\Gh}
 \lesssim\theta(\kappa/h)^{1/2}\norme v_h$.
Since $\beta\rho s=\theta\rho$, $\beta s=\theta$, and
$\theta\le\sqrt\theta$, these estimates bound the part containing $r_T$
by
\[
 C\delh^{m+1}C_u(\kappa\theta/h)^{1/2}\norme v_h .
\]
For the flux residual, the continuous trace estimate and
\reff{eq:poincare} give
$\norm v_{\Gh}\le C_\Lambda\kappa^{-1/2}\norme v_h$.
Combining this with $s=(1-\theta)h/(\gamma\kappa)$ yields
\[
 \rho s\norm v_{\Gh}
 +s\norm{\kappa\p_{\nh}v}_{\Gh}
 \lesssim(1-\theta)\kappa^{-1/2}
 (C_\Lambda+h^{1/2})\norme v_h .
\]
Multiplication by $\norm{r_\Lambda}_{\Gh}$ proves
\reff{eq:consbound}.
Finally,
\[
 \frac\theta h=\frac\beta{\beta h+\gamma\kappa}
 \le\frac\beta{\gamma\kappa}
 \qquad (\beta<\infty).
\]
This proves the estimate for fixed finite $\beta$, with constants that
may depend on $\beta$ and $\Lambda$.
\Endproof

\subsection{The energy error.}

\begin{thm}[Energy error estimate]
\label{thm:err}
Under the hypotheses of \Thm{\ref{thm:stab}} and \Asm{\ref{asm:ext}},
let $m\le p-1$ and $h\le h_0$.  Then
\be
  \label{eq:err}
  \norme{u-u_h}_h\lesssim h^p\norm u_{H^{p+1}}
  +\delh^{m+1}C_u\sqrt\kappa
   \left[\left(\frac\theta h\right)^{1/2}
        +(1-\theta)(C_\Lambda+h^{1/2})\right] .
\ee
Constants are independent of $h$, $\beta$, and the position of the cut
when the Poincar\'e constant $\Lambda$ is uniformly bounded.
For fixed $0\le\beta<\infty$ this implies
\be
  \label{eq:errfinite}
  \norme{u-u_h}_h
    \lesssim_\beta h^p\norm u_{H^{p+1}}+\delh^{m+1}C_u .
\ee
Thus $2m+2\ge p$ gives the optimal energy order for fixed finite $\beta$,
while $2m+\tfrac32\ge p$ remains sufficient uniformly up to the Dirichlet
limit.
\end{thm}

\Proof
Let $I_hu$ be an interpolant on the active mesh and set
$\eta=u-I_hu$ and $e_h=I_hu-u_h$.
Standard interpolation and trace estimates on the full elements give
\[
 \norm{D^k\eta}_{\Gh}
  \lesssim h^{p+1/2-k}\norm u_{H^{p+1}}
  \qquad (0\le k\le p).
\]
They also give $\norme\eta_h\lesssim h^p\norm u_{H^{p+1}}$, with the
extension used in all these norms.  For a zero-mean space, we adjust the
interpolant by a constant without changing the approximation order.

We split $A_h=A_h^0+B_h$ as in \reff{eq:split} and first consider
$A_h^0(\eta,v)$.  The volume and ghost-penalty terms have the usual
interpolation bounds, while the boundary terms are controlled by
\[
 \theta\norm v_{\Gh}
 \le(\theta/\rho)^{1/2}\norme v_h,
 \qquad
 \norm{\kappa\p_{\nh}v}_{\Gh}
 \lesssim(\kappa/h)^{1/2}\norme v_h,
 \qquad s\le1/\rho .
\]
The first estimate is interpreted as zero at $\theta=0$.
Together with the interpolation bounds, these estimates give
$\abs{A_h^0(\eta,v)}\lesssim h^p\norm u_{H^{p+1}}\norme v_h$.
In particular, no inverse power of $\theta$ is needed.

To estimate $B_h$, write
$D_T(\eta)=\Tm(\eta)-\eta$ and
$D_\Lambda(\eta)=\Lm(\eta)-\kappa\p_{\nh}\eta$.
Since $m\le p-1$ and $\delh\lesssim h^2$, the same trace estimates imply
\[
 \begin{aligned}
 \norm{D_T(\eta)}_{\Gh}
  &\lesssim \delh h^{p-1/2}\norm u_{H^{p+1}},\\
 \norm{D_\Lambda(\eta)}_{\Gh}
  &\lesssim\kappa
    \left(\thh h^{p-1/2}+\delh h^{p-3/2}\right)
       \norm u_{H^{p+1}}
   \lesssim\kappa h^{p+1/2}\norm u_{H^{p+1}} .
 \end{aligned}
\]
For $m=0$, the first expression vanishes and the $\delh$ term in the
second expression is absent.  We apply the continuous trace estimate to
$v$ in the first term of \reff{eq:Bh} and use the boundary penalty and
the discrete flux estimate in the other three terms.  This gives
\[
 \abs{B_h(\eta,v)}
 \lesssim\left(h^{p+1/2}+h^{p+1}\right)
       \norm u_{H^{p+1}}\norme v_h .
\]
Its constant depends on the stated geometric and Poincar\'e constants.
Consequently,
$\abs{A_h(\eta,v)}\lesssim h^p\norm u_{H^{p+1}}\norme v_h$.

Coercivity and the discrete equation now give
\[
 \norme{e_h}_h^2\lesssim
  \abs{A_h(\eta,e_h)}
  +\abs{A_h(u,e_h)-L_h(e_h)} .
\]
Applying \Cor{\ref{cor:cons}} and then the triangle inequality proves
\reff{eq:err}.
The last assertions follow from $\delh\lesssim h^2$ and
$\theta/h\le\beta/(\gamma\kappa)$ for finite $\beta$.
\Endproof

\subsection{The $L^2$ error.}

\begin{rem}[Expected $L^2$ estimate]
\label{pro:l2}
The expected estimate is
\be
  \label{eq:l2}
  \norm{u-u_h}_{\Omh}\lesssim h^{p+1}\norm{u}_{H^{p+1}}
       +\delh^{m+1}C_u .
\ee
For the Dirichlet problem, the corresponding estimate is proved in
\cite[Thm.~3.2]{BHL18} under the assumptions made there.
We do not prove \reff{eq:l2} for finite $\beta$ in this paper.
Since the trial-side corrections make the form non-symmetric, a duality
proof requires a separate adjoint-consistency estimate.
It must also control the geometric remainders using the regularity
available for the dual solution, usually $H^2$.
The $W^{m+2,\infty}$ remainder bound above does not provide this directly.

These points remain open after the improvement of the energy estimate.
For example, $p=2$ and $m=0$ give order two in energy at fixed finite
$\beta$, but the geometric term in \reff{eq:l2} would still be of order
two.  An additional Taylor term is therefore needed to obtain the
expected third order in $L^2$, in agreement with the numerical results.
\end{rem}

\subsection{A symmetric preconditioner.}
\label{ss:precond}

Setting $m=0$ and replacing $\nst$ by $\nh$ in \reff{eq:Ah} gives the
symmetric form
\be
  \label{eq:A0}
  \begin{split}
  A_h^0(u,v):={}&(\kappa\nabla u,\nabla v)_{\Omh}
       +(cu,v)_{\Omh}+\jh(u,v)\\
   &-\theta\big(\pd{\kappa\p_{\nh}u,v}_{\Gh}
                   +\pd{u,\kappa\p_{\nh}v}_{\Gh}\big)\\
   &+\theta\rho\pd{u,v}_{\Gh}
       -s\pd{\kappa\p_{\nh}u,\kappa\p_{\nh}v}_{\Gh}.
  \end{split}
\ee
This is the uncorrected Juntunen--Stenberg form on $\Omh$, with the
ghost penalty added.
Define
\be
  \label{eq:split}
  A_h=A_h^0+B_h,
  \qquad D_T(u):=\Tm(u)-u,
  \qquad D_\Lambda(u):=\Lm(u)-\kappa\p_{\nh}u .
\ee
Then
\be
  \label{eq:Bh}
  \begin{split}
  B_h(u,v)={}&\rho s\pd{D_\Lambda(u),v}_{\Gh}
       +\theta\rho\pd{D_T(u),v}_{\Gh}\\
       &-\theta\pd{D_T(u),\kappa\p_{\nh}v}_{\Gh}
       -s\pd{D_\Lambda(u),\kappa\p_{\nh}v}_{\Gh}.
  \end{split}
\ee

\begin{pro}[A small perturbation of a symmetric form]
\label{pro:pert}
Let $\gamma>4C_I$ and $0\le m\le p-1$.
Under the mesh assumptions, the symmetric form satisfies
\be
  \label{eq:A0coer}
  A_h^0(v,v)\ge C_0\norme{v}_h^2,
  \qquad
  \abs{A_h^0(u,v)}\le C_1\norme{u}_h\norme{v}_h .
\ee
Under \reff{eq:poincare},
\be
  \label{eq:Bbound}
  \abs{B_h(u,v)}\le C_Bc_\Gamma h^{1/2}
      \norme{u}_h\norme{v}_h .
\ee
The constants are independent of the cut position, although $C_B$ may
depend on $\Lambda$.
\end{pro}

\Proof
Coercivity follows from the argument preceding ${\cal R}$ in
\reff{eq:stabsplit}, and continuity follows from \reff{eq:disc} and the bounds
$\theta\le1$ and $s\le1/\rho$.
By \reff{eq:tpert} and \reff{eq:E},
\[
 \norm{D_T(u)}_{\Gh}\lesssim\frac{\delh}{\sqrt{\kappa h}}\norme{u}_h,
 \qquad
 \norm{D_\Lambda(u)}_{\Gh}
 \lesssim\Big(\thh+\frac{\delh}{h}\Big)
         \Big(\frac{\kappa}{h}\Big)^{1/2}\norme{u}_h.
\]
Substitution in \reff{eq:Bh} bounds the two trace-correction terms by
$C(\delh/h)\norme{u}_h\norme{v}_h$ and the flux-flux term by
$C(\thh+\delh/h)\norme{u}_h\norme{v}_h$.
The first term is bounded by
$CC_\Lambda(\thh+\delh/h)h^{-1/2}\norme{u}_h\norme{v}_h$
using \reff{eq:energytrace}.
The geometric estimates now give \reff{eq:Bbound}. \Endproof

\begin{cor}[Preconditioning by the uncorrected form]
\label{cor:precond}
Let $q(h):=C_Bc_\Gamma h^{1/2}/C_0<1$.
In the inner product induced by $A_h^0$,
\be
  \label{eq:iter}
  \norm{I-(A_h^0)^{-1}A_h}_{A_h^0}\le q(h).
\ee
Thus preconditioned Richardson iteration converges with contraction
factor at most $q(h)$, and the operator condition number in this norm is
at most $(1+q)/(1-q)$.
\end{cor}

\Proof
For a discrete function $x$,
\[
 \norm{(A_h^0)^{-1}B_hx}_{A_h^0}
 =\sup_{y\ne0}\frac{\abs{B_h(x,y)}}{\norm{y}_{A_h^0}}
 \le\frac{C_Bc_\Gamma}{C_0}h^{1/2}\norm{x}_{A_h^0}.
\]
This proves the operator bound.  When the preconditioned operator acts
on a unit vector, the resulting norm lies between $1-q$ and $1+q$, which
gives the condition-number estimate. \Endproof

The symmetric matrix uses the same quadrature points as the corrected
method but does not require $\rh$ or $\nst$.
The correction can then be applied through matrix--vector products with
$B_h$, as illustrated by the numerical results in \S\ref{ss:prec}.
An analogous splitting is available for the interface form, subject to
the phasewise Poincar\'e bounds used in \S\ref{sec:iface}.

\Section{The interface problem with imperfect bonding.}
\label{sec:iface}

\subsection{The problem.}

Let $\Gamma$ be a closed $C^2$ interface in a fixed domain $\Omega$,
dividing $\Omega$ into an enclosed subdomain $\Omega_1$ and a surrounding
subdomain $\Omega_2$.  We take the normal $\V n$ to point from $\Omega_1$
into $\Omega_2$ and let $\kappa_1,\kappa_2>0$.
For a function $v$ defined on both sides, set $\jm{v}:=v_1-v_2$.
Given weights $w_1+w_2=1$, we also define
$\av{v}:=w_1v_1+w_2v_2$ and the complementary average
$\av{v}':=w_2v_1+w_1v_2$.  With this notation, we consider the problem
\be
  \label{eq:imodel}
  \begin{array}{rcll}
   -\nabla\cdot(\kappa_i\nabla u)&=&f&\hbox{in }\Omega_i,\ i=1,2,\\[2pt]
   \jm{\kappa\p_{\V n}u}&=&g_N&\hbox{on }\Gamma,\\[2pt]
   \beta_\Gamma\jm{u}+\av{\kappa\p_{\V n}u}&=&g_\Gamma&\hbox{on }\Gamma ,
  \end{array}
\ee
with a Dirichlet condition on the fitted outer boundary $\p\Omega$.
The third equation describes imperfect bonding through a spring layer
with stiffness $\beta_\Gamma\in[0,\infty]$.
For homogeneous data, the flux is proportional to the jump in $u$.
The limit $\beta_\Gamma=\infty$ gives perfect bonding, with a prescribed
jump $\jm{u}=g_D$ obtained by scaling the data as in the boundary Robin
problem.  At $\beta_\Gamma=0$, the two fluxes are prescribed independently
and vanish when $g_N=g_\Gamma=0$.
This is the scalar counterpart of the spring-layer law of Zhong and Meguid
\cite{ZM99}, which A.~Hansbo and P.~Hansbo \cite{HH04} analysed and
discretised in elasticity.  We return to this relation in \S\ref{sec:HH}.

The bonding law has the form of the Robin condition, with $\jm{u}$ replacing
the boundary value and $\av{\kappa\p_{\V n}u}$ replacing the boundary flux.
We can therefore apply \Lem{\ref{lem:subst}} directly, although the
flux-jump condition requires an additional correction.

\subsection{Three transports instead of two.}

We construct $\Gh$ as in \S\ref{sec:geom}, so that it divides $\Omega$
into two polytopal subdomains $\Oh1$ and $\Oh2$.
Let $\Th^\Gamma$ denote the cut elements.
We use the Hansbo--Hansbo space, with separate unknowns on the two sides,
\be
  \label{eq:ispace}
  \Vh:=\Vh^1\times\Vh^2,\qquad
  \Vh^i:=\set{v\in C^0:\ v|_T\in\mP_p(T)\ \forall T\in\Th^{{\rm act},i}} ,
\ee
where each factor is defined on the elements meeting its subdomain.
Thus each cut element carries two independent sets of degrees of freedom.
We choose the harmonic weights and mean coefficient
\be
  \label{eq:weights}
  w_1=\frac{\kappa_2}{\kappa_1+\kappa_2},\quad
  w_2=\frac{\kappa_1}{\kappa_1+\kappa_2},\quad
  \kbar=\frac{2\kappa_1\kappa_2}{\kappa_1+\kappa_2},\quad
  \rho=\frac{\gamma\kbar}{h},
\ee
following Dryja \cite{D03} and Burman and Zunino \cite[(2.22)]{BZ06},
with a related discussion in \cite[\S2]{BGSS18}.
With $\omega:=w_1\kappa_1=w_2\kappa_2=\kbar/2$, we have
$\av{\kappa\nabla v\cdot\V n}
 =\omega(\nabla v_1+\nabla v_2)\cdot\V n$.
Thus the weighted flux average is $\kbar$ times the arithmetic average
of the normal gradients \cite[\S2.2]{BZ06}.
Although other weights with $w_1+w_2=1$ can be used in the analysis
\cite[Rem.~2.13]{BZ06}, the present choice controls the dependence on the
coefficient contrast.  This property is used in \Lem{\ref{lem:imis}}.
Integration by parts on the polygonal subdomains gives
\be
  \label{eq:iibp}
  \sum_i(\kappa_i\nabla u,\nabla v)_{\Oh i}
   -\pd{\av{\kappa\p_{\nh}u},\jm{v}}_{\Gh}
   -\pd{\jm{\kappa\p_{\nh}u},\av{v}'}_{\Gh}=(f,v)_{\Omega} ,
\ee
where we used the identity
$a_1b_1-a_2b_2=\av{a}\jm{b}+\jm{a}\av{b}'$.
In addition to $\Tm$ and $\Lm$, we need a transport of the flux jump,
\be
  \label{eq:Pim}
  \Pim(w)(x):=\Tm\big(\jm{\kappa\nabla w}\big)(x)\cdot\nst(x) ,
\ee
which is the analogue of \reff{eq:Lm} for the flux jump.
We expand $\jm{\kappa\nabla w}$ componentwise along $\nh$ and contract the
result with the fixed vector $\nst(x)$, so that only $w$ is differentiated.
Applying \Lem{\ref{lem:rem}} on each side then shows that $\Pim(u)$
approximates $\jm{\kappa\p_{\V n}u}(\ph(x))$ with an error of order
$\delh^{m+1}$.

\subsection{The method.}

As in the boundary problem, we introduce the weights
\be
  \label{eq:ithetas}
  \thg:=\frac{\beta_\Gamma}{\beta_\Gamma+\rho}\in[0,1],
  \qquad
  \sg:=\frac{1}{\beta_\Gamma+\rho},
  \qquad
  \rho=\frac{\gamma\kbar}{h} ,
\ee
which satisfy $\thg=\beta_\Gamma\sg$ and $1-\thg=\rho\sg$.
We can therefore apply \Lem{\ref{lem:subst}} with $U=\jm{u}$,
$\Lambda=\av{\kappa\p_{\V n}u}$, $G=g_\Gamma$ and $\beta=\beta_\Gamma$.
The discrete problem is to find $u_h\in\Vh$ such that
$A_h^\Gamma(u_h,v)=L_h^\Gamma(v)$ for all $v\in\Vh$, where
\be
  \label{eq:iAh}
  \begin{array}{rcl}
  A_h^\Gamma(u,v)&:=&\displaystyle\sum_i(\kappa_i\nabla u,\nabla v)_{\Oh i}
      +\jh(u,v)
      -\pd{\av{\kappa\p_{\nh}u},\jm{v}}_{\Gh}\\[8pt]
   &&{}-\ \pd{\jm{\kappa\p_{\nh}u}-\Pim(u),\ \av{v}'}_{\Gh}\\[6pt]
   &&{}+\ \rho\sg\pd{\Lm(u),\jm{v}}_{\Gh}
       +\thg\rho\pd{\jm{\Tm(u)},\jm{v}}_{\Gh}\\[6pt]
   &&{}-\ \thg\pd{\jm{\Tm(u)},\av{\kappa\p_{\nh}v}}_{\Gh}
       -\sg\pd{\Lm(u),\av{\kappa\p_{\nh}v}}_{\Gh},
  \end{array}
\ee
\be
  \label{eq:iLh}
  L_h^\Gamma(v):=(f,v)_{\Omega}
    +\pd{g_N\circ\ph,\av{v}'}_{\Gh}
    +\rho\sg\pd{g_\Gamma\circ\ph,\jm{v}}_{\Gh}
    -\sg\pd{g_\Gamma\circ\ph,\av{\kappa\p_{\nh}v}}_{\Gh} ,
\ee
Here $\Lm(u):=\Tm(\av{\kappa\,\nst\!\cdot\!\nabla u})$ transports
the weighted average flux, while the ghost penalty $\jh$ acts separately
on the two active meshes.  On each side, it is applied to faces belonging
to cut elements.

The second line of \reff{eq:iAh} corrects the flux jump, while the remaining
terms follow from the Robin method in \S\ref{sec:method}.
The value of $\beta_\Gamma$ enters only through $\thg$ and $\sg$.
At $\beta_\Gamma=\infty$, we have $\thg=1$ and $\sg=\rho\sg=0$, so that
the two $\Lm$ terms vanish.
With the corresponding scaling of the data, this gives the
Hansbo--Hansbo method \cite{HH02} for $\jm{u}=g_D$, with the jump corrected
by $\Tm$ and the flux jump by $\Pim$.
At $\beta_\Gamma=0$, the values $\thg=0$, $\sg=1/\rho$ and $\rho\sg=1$
remove the $\Tm$ terms, and the two flux conditions prescribe the
one-sided fluxes.  For $0<\beta_\Gamma<\infty$, both transports contribute
with weights $\thg$ and $1-\thg$.
The experiments in \S\ref{ss:imis} cover nine orders of magnitude in
$\beta_\Gamma$.

\begin{rem}[The flux-jump term is the misalignment correction]
\label{rem:keep}
On the exact interface, the term
$\pd{\jm{\kappa\p_{\V n}u},\av{v}'}_\Gamma$ vanishes when the flux is
continuous.  Its counterpart on $\Gh$ must nevertheless be retained,
as shown by the exact identity on $\Gh$,
\be
  \label{eq:tangential}
  \jm{\kappa\p_{\nh}u}=\Pi_0(u)+\jm{\kappa\nabla u}\cdot(\nh-\nst).
\ee
Here $\Pi_0(u)$ approximates $g_N\circ\ph$ with an error of order
$\delh$, while the second term contains the tangential flux jump and a
normal mismatch of size $\thh=O(h)$.
The tangential flux jump is generally non-zero, since on $\Gamma$ it
satisfies
$\jm{\kappa\nabla u}\cdot\V t
 = (\kappa_1-\kappa_2)\p_{\V t}u_1
   +\kappa_2\p_{\V t}\jm{u}$.
This term therefore depends on the solution and cannot be replaced by
data.  A fixed subdivision of each chord reduces its size but does not
improve its order in $h$.
We instead retain the computable flux jump $\jm{\kappa\p_{\nh}u}$ in the
method.  After integration by parts, the correction by $\Pim(u)$ leaves
only the residual $\Pim(u)-g_N\circ\ph$, which is of order
$\delh^{m+1}$ at the exact solution.
\end{rem}

\subsection{Consistency and error estimate.}

In this section we use the interface energy norm
\be
 \label{eq:inorm}
 \norme{v}_h^2:=\sum_{i=1}^2\kappa_i\norm{\nabla v_i}^2_{\Oh i}
    +\jh(v,v)+\thg\rho\norm{\jm{v}}^2_{\Gh}.
\ee
Here $\jh=\sum_i j_{h,i}$ is the sum of the two ghost penalties.
At $\beta_\Gamma=0$, we fix the mean on the enclosed phase by imposing
zero mean there on both the trial and test spaces.
For the error estimate, we use the same normalization for the exact
solution.

We define the two transported residuals
\be
  \label{eq:iR}
  R_\Gamma(u):=\beta_\Gamma\jm{\Tm(u)}+\Lm(u)-g_\Gamma\circ\ph,
  \qquad
  R_N(u):=\Pim(u)-g_N\circ\ph .
\ee

\begin{thm}[Consistency]
\label{thm:icons}
Let $0\le\beta_\Gamma<\infty$ and let $u$ solve \reff{eq:imodel}, with each
$u_i$ satisfying \Asm{\ref{asm:ext}} on its side.
Then, for all $v\in\Vh$,
\be
  \label{eq:icons}
  A_h^\Gamma(u,v)-L_h^\Gamma(v)
   =\rho\sg\pd{R_\Gamma(u),\jm{v}}_{\Gh}
    -\sg\pd{R_\Gamma(u),\av{\kappa\p_{\nh}v}}_{\Gh}
    +\pd{R_N(u),\av{v}'}_{\Gh} ,
\ee
and $\norm{R_\Gamma(u)}_{\Gh}+\norm{R_N(u)}_{\Gh}\lesssim
 \delh^{m+1}(\beta_\Gamma+\kappa_1+\kappa_2)\,C_u$.
\end{thm}

\Proof
As in the proof of \Thm{\ref{thm:cons}}, we use \reff{eq:iibp} to cancel
the volume terms and the two boundary terms.
The flux-jump correction and its datum then leave
$\pd{R_N(u),\av{v}'}$.
Applying \Lem{\ref{lem:subst}} to the remaining terms, with
$U=\jm{u}$ and $\Lambda=\av{\kappa\p_{\V n}u}$, gives the first two terms
of \reff{eq:icons}.  The residual bound follows from
\Lem{\ref{lem:rem}} on each side.
\Endproof

\begin{thm}[Energy error estimate]
\label{thm:ierr}
Assume the phasewise Poincar\'e bounds of \Lem{\ref{lem:imis}} and
\Asm{\ref{asm:ext}} on each side, and let $0\le m\le p-1$.
For $\gamma$ sufficiently large and $h\le h_0$, the discrete problem has a
unique solution when $\beta_\Gamma>0$.
When $\beta_\Gamma=0$, impose the Neumann compatibility condition and fix
the mean in the enclosed phase, restricting both trial and test spaces
as specified above.  Then
\be
  \label{eq:ierr}
  \Big(\sum_i\norm{\sqrt{\kappa_i}\nabla(u-u_h)}^2_{\Oh i}
       +\thg\rho\norm{\jm{u-u_h}}^2_{\Gh}\Big)^{1/2}
  \ \lesssim\ h^p\Big(\sum_i\kappa_i\norm{u_i}^2_{H^{p+1}(U_i)}\Big)^{1/2}
   +h^{-1/2}\delh^{m+1}C_u ,
\ee
Here $U_i$ is the fixed extension neighbourhood of $\overline{\Omega_i}$,
and $C_u$ carries the same coefficient weights.
The constants and $h_0$ may depend on the phasewise Poincar\'e bounds but
are independent of the position of $\Gamma$ relative to the mesh.
They are also independent of $\beta_\Gamma$ and the contrast
$\kappa_2/\kappa_1$, provided that the Poincar\'e bounds are uniform in
these parameters.
\end{thm}

\begin{rem}[What contrast independence means here]
\label{rem:contrast}
Since the norms on both sides of \reff{eq:ierr} contain coefficient
weights, contrast independence refers to the constants in these norms,
subject to the stated Poincar\'e bounds.
For perfect bonding, Burman, Guzm\'an, S\'anchez and Sarkis
\cite[Thm.~4.3]{BGSS18} prove a stronger flux estimate whose constant is
independent of the coefficients.
Their proof uses a discrete extension from one phase to the whole domain.
We have not examined whether boundary value correction preserves that
estimate.
\end{rem}

\Sketch
We follow the coercivity proof of \Thm{\ref{thm:stab}}, using the weights
\reff{eq:weights} to control the coefficient factors \cite{HH02,BZ06}.
The last two lines of \reff{eq:iAh} are treated as in the boundary problem,
with $\jm{v}$ in place of $v$.
The second line gives the additional misalignment term
$-\pd{\jm{\kappa\p_{\nh}v}-\Pi_0(v),\av{v}'}=
 -\pd{\jm{\kappa\nabla v}\cdot(\nh-\nst),\av{v}'}$, which is bounded in \Lem{\ref{lem:imis}} below.
Consistency follows from \Thm{\ref{thm:icons}}, and we bound the first
two terms of \reff{eq:icons} as in \Cor{\ref{cor:cons}}.
For the flux-jump residual, we retain the coefficient weights before
applying the Cauchy--Schwarz inequality.
A Taylor remainder from phase $i$ paired with the trace from phase $j$
has coefficient $\kappa_i w_{3-j}\kappa_j^{-1/2}$.
Extracting the solution weight $\kappa_i^{1/2}$ leaves the factor
\[
 \kappa_i^{1/2}w_{3-j}\kappa_j^{-1/2}\le1,
 \qquad i,j\in\{1,2\}.
\]
The calculation is the same as in \Lem{\ref{lem:imis}} below, so the
continuous trace estimate gives an $O(\delh^{m+1})$ bound for this residual
in the weighted norms.
We treat the average-flux residual coefficient by coefficient as well,
and use the boundary penalty for the trace residual.
Together, these estimates imply the stated uniform
$O(h^{-1/2}\delh^{m+1})$ bound.
The perfect-bonding case follows by taking the scaled-data limit, while
the interpolation argument is the phasewise version of
\Thm{\ref{thm:err}}.
\Endproof

\begin{lem}[The interface misalignment remainder]
\label{lem:imis}
Let the weights be those of \reff{eq:weights} and assume, as in
\reff{eq:poincare}, that $\kappa_i\norm{v_i}^2_{\Oh i}\le\Lambda^2\norme{v}_h^2$
for $i=1,2$.  Then
$$
  \abs{\pd{\jm{\kappa\nabla v}\cdot(\nh-\nst),\av{v}'}_{\Gh}}
   \ \le\ C\,c_\Gamma\,h^{1/2}\,\norme{v}_h^2 ,
$$
where $C$ depends on $\Lambda$ but, for a fixed bound on $\Lambda$, is
independent of $h$, $\beta_\Gamma$ and the contrast $\kappa_2/\kappa_1$.
\end{lem}

\Proof
As in \Lem{\ref{lem:mis}}, we apply the discrete trace inequality to the
gradients and the continuous trace inequality \reff{eq:ctrace} to the
functions.  This gives
$$
  \begin{array}{l}
  \displaystyle
  \norm{\jm{\kappa\nabla v}}_{\Gh}\le\sum_i\kappa_i\norm{\nabla v_i}_{\Gh}
   \lesssim h^{-1/2}\sum_i\kappa_i^{1/2} E_i,
  \\[10pt]\displaystyle
  \norm{\av{v}'}_{\Gh}\le\sum_i w_{3-i}\norm{v_i}_{L^2(\Gh)}
   \lesssim\sum_i w_{3-i}\,\kappa_i^{-1/2}\Lambda_i ,
  \end{array}
$$
where
$E_i:=(\kappa_i\norm{\nabla v_i}^2_{\Oh i}+j_{h,i}(v_i,v_i))^{1/2}
\le\norme{v}_h$ and
$\Lambda_i:=(1+\Lambda^2)^{1/2}\norme{v}_h$.  The product contains the
four coefficients $\kappa_i^{1/2}w_{3-j}\kappa_j^{-1/2}$.  The two diagonal
ones are $w_2$ and $w_1$, both at most $1$.  The two off-diagonal ones are
both equal to
$$
  \Big(\frac{\kappa_1}{\kappa_2}\Big)^{1/2}\frac{\kappa_2}{\kappa_1+\kappa_2}
  =\frac{\sqrt{\kappa_1\kappa_2}}{\kappa_1+\kappa_2}\ \le\ \frac12
$$
by the arithmetic--geometric mean inequality.
Thus all four coefficients are bounded by $1$ independently of the
contrast, and multiplication by $\thh\le2c_\Gamma h$ gives the result.
The harmonic weights \cite{BZ06} are essential here, since arithmetic
weights would introduce the off-diagonal factors
$\frac12(\kappa_1/\kappa_2)^{1/2}$ and its reciprocal.
\Endproof

\begin{rem}[The debonded limit]
\label{rem:ismall}
At $\beta_\Gamma=0$, the enclosed phase solves a pure Neumann problem,
whose solution is determined only up to a constant.
We therefore fix its mean before applying Poincar\'e's inequality, as in
the related normalisation used in \cite[Thm.~4.3]{BGSS18}.
For $\beta_\Gamma>0$, the bonding term controls this constant, but that
control weakens as $\beta_\Gamma$ tends to zero.
Thus \Thm{\ref{thm:ierr}} does not assert uniformity in this limit unless
the phasewise Poincar\'e constants remain bounded.
\end{rem}

\Section{Imperfect bonding as the interface Robin condition.}
\label{sec:HH}

A.~Hansbo and P.~Hansbo \cite{HH04} considered elasticity problems with
strong and weak discontinuities, coupling the two phases through a
cohesive law that includes the limit of perfect bonding.
The third equation in \reff{eq:imodel} is the scalar analogue of this law.
Our formulation \reff{eq:iAh} uses the same doubled unknowns and weighted
averages, with Taylor transports and a flux-jump correction to account for
the surrogate interface.

Boundary value correction also has a direct interpretation as an interface
compliance.  To see this, take $m=1$ and assume continuous flux, so that
$g_N=0$.  Then, on $\Gh$,
$$
  \jm{{\cal T}_1(u)}=\jm{u}+\rh\,\jm{\p_{\nh}u}
   =\jm{u}+{\cal K}\,\kappa_1\p_{\nh}u_1+O((\thh+\delh)\delh) ,
$$
where ${\cal K}:=\rh(\kappa_1^{-1}-\kappa_2^{-1})$.
Thus the transported condition $\jm{{\cal T}_1(u)}=g_D\circ\ph$ gives,
to this order,
\be
  \label{eq:compliance}
  \jm{u}=g_D\circ\ph-{\cal K}\,\kappa_1\p_{\nh}u_1 .
\ee
This has the algebraic form of an imperfect bonding law with compliance
${\cal K}$, although the geometric compliance can have either sign.
The geometric offset thus acts as an artificial compliance, weighted by
the jump of the inverse diffusivity.
Since $\rh=O(h^2)$, this contribution has order $h^2$ and, at this order,
adds to the physical compliance $1/\beta_\Gamma$.

Two features distinguish this correction from a physical compliance.
Because ${\cal K}$ may change sign with $\rh$, we treat it as a
perturbation in the stability analysis.
Moreover, ${\cal K}=0$ when $\kappa_1=\kappa_2$, so the first-order
correction vanishes under the assumed flux continuity.
At a boundary, the weights $\theta$ and $1-\theta$ give a corresponding
interpretation.  The Dirichlet limit corrects the value and the Neumann
limit corrects the flux, while both corrections contribute for a Robin
condition.

\Section{Numerical examples.}
\label{sec:num}

\begin{figure}[!htbp]\centering
\includegraphics[width=0.46\textwidth]{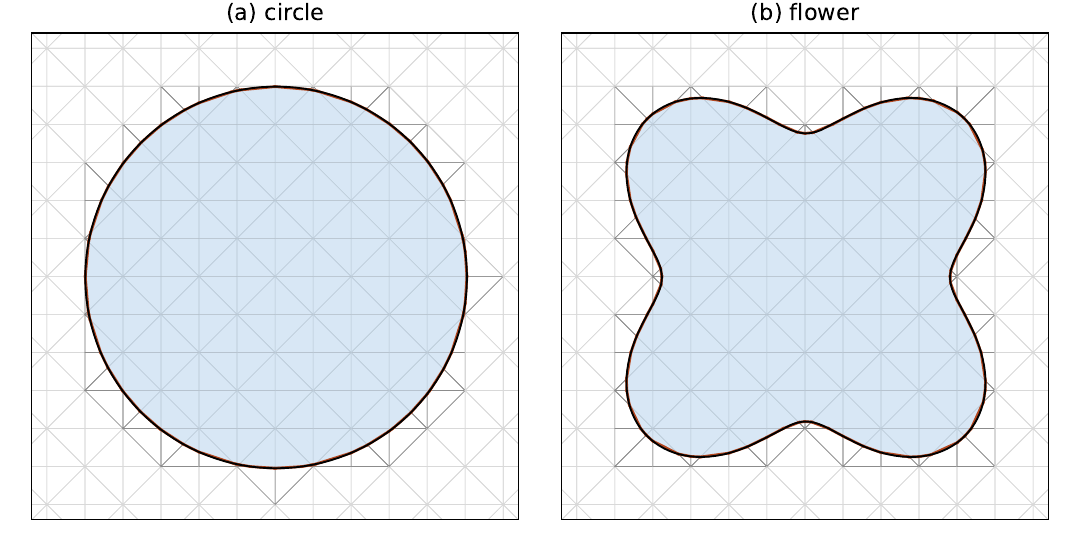}
\caption{The two geometries at $N=16$: the background mesh, the active elements (shaded), the surrogate boundary $\Gh$ (thick) and the exact boundary $\Gamma$ (thin).}
\label{fig:mesh}
\end{figure}

We use structured triangulations of the unit square and continuous
$\mP_p$ elements, with one chord per cut element for the surrogate boundary.
The ghost penalty \reff{eq:ghost} uses $\gamma_g=5\cdot10^{-2}$, and the
Nitsche parameter is set to $\gamma=20$.
The first geometry is a circle of radius $0.3137$, centred at
$(0.5017,0.4983)$.
The second is a four-petal flower with level set
$\Phi=r^2-R^2+c\cos(4\varphi)$, where $R=0.31$ and $c=0.04$, about the same
centre.  Neither boundary is aligned with the mesh, as shown in
\Fig{\ref{fig:mesh}}.

For the boundary value problems, we take
$u=\cos(2x)\sin(2y)$ with $\kappa=1$ and $c=1$, compute $f$ and $g$ from
this solution, and evaluate the boundary datum at $\ph(x)\in\Gamma$.
Cut cells are integrated using a centroid fan with signed triangle areas.
With sufficiently accurate triangle quadrature, this gives exact
integration of polynomials on a simple polygon.
On each element, we obtain the nodal basis and the derivatives required by
the Taylor operators by inverting the Vandermonde matrix of a scaled
monomial basis.

We first apply a patch test on a curved boundary.
For a polynomial solution of degree $p$ and $m=p$, the operators
$\Tm$, $\Lm$ and $\Pim$ are exact, so the residuals \reff{eq:R} and
\reff{eq:iR} vanish.
This algebraic test uses $m=p$, whereas the error analysis assumes
$m\le p-1$.
The relative $H^1$ error is below $2\cdot10^{-15}$ for $p=1$,
$4\cdot10^{-13}$ for $p=2$ and $6\cdot10^{-10}$ for $p=3$.
These bounds hold for $\beta\in\set{0,1,10^3,\infty}$ in both the boundary
and interface tests, with the larger error for $p=3$ consistent with the
conditioning of the monomial Vandermonde matrix.
The patch test thus checks the geometry maps, Taylor operators, Nitsche
terms, ghost penalty and quadrature together.

We also measure the geometric errors for $N=12,\dots,192$.
The ratio $\delh/h^2$ stays below $0.26$ on the circle and $1.02$ on the
flower, while the corresponding bounds for $\thh/h$ are $1.18$ and $5.78$.
These measurements are consistent with the orders in
\Lem{\ref{lem:geom}}.

\subsection{The Robin family: convergence and the Taylor order.}
\label{ss:conv}

\Tab{\ref{tab:conv}} reports least-squares convergence rates in the $H^1$
seminorm and in $L^2$ over the four meshes $N=12,24,48,96$, for $p=1,2,3$,
$m=0,1,2$ and the three representative values $\beta=0$ (Neumann),
$\beta=1$ (Robin) and $\beta=\infty$ (Dirichlet).

\begin{table}[!htbp]\centering
\begin{tabular}{llrrrrrr}\toprule
& & \multicolumn{2}{c}{$\beta=0$} & \multicolumn{2}{c}{$\beta=1$} & \multicolumn{2}{c}{$\beta=\infty$} \\
$p$ & $m$ & $H^1$ & $L^2$ & $H^1$ & $L^2$ & $H^1$ & $L^2$ \\ \midrule
1 & 0 & 0.98 & 1.79 & 0.98 & 1.92 & 0.99 & 2.01 \\
1 & 1 & 0.98 & 1.79 & 0.98 & 1.92 & 0.99 & 2.02 \\
\midrule
2 & 0 & 2.15 & 1.77 & 2.20 & 1.82 & 1.54 & 1.87 \\
2 & 1 & 2.22 & 2.96 & 2.22 & 3.16 & 2.22 & 3.29 \\
2 & 2 & 2.22 & 2.96 & 2.22 & 3.16 & 2.22 & 3.29 \\
\midrule
3 & 0 & 1.81 & 1.78 & 1.85 & 1.82 & 1.33 & 1.90 \\
3 & 1 & 3.40 & 4.07 & 3.40 & 4.24 & 3.41 & 4.06 \\
3 & 2 & 3.40 & 4.00 & 3.40 & 4.19 & 3.41 & 4.36 \\
\bottomrule
\end{tabular}
\caption{Least-squares convergence rates on the circle over $N=12,\dots,96$.  Optimal is $p$ in $H^1$ and $p+1$ in $L^2$.}
\label{tab:conv}
\end{table}

\begin{figure}[!htbp]\centering
\includegraphics[width=0.80\textwidth]{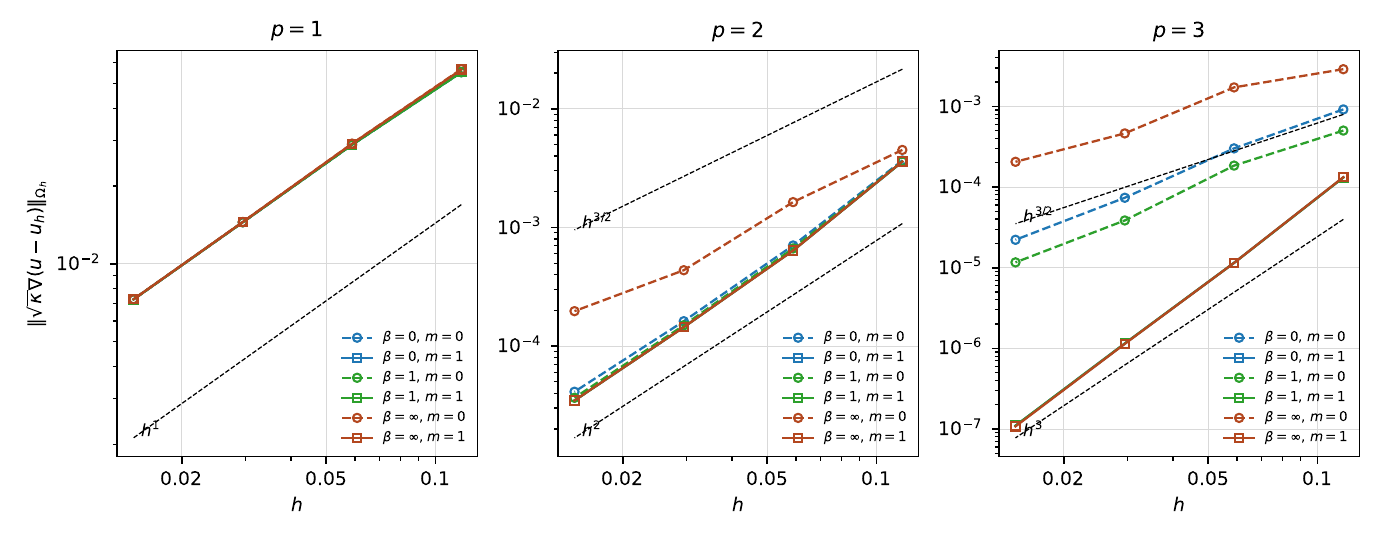}
\caption{Energy-norm convergence on the circle for $p=1,2,3$.  Dashed curves use $m=0$ and solid curves use $m=1$.  The reference slopes are $h^p$ and $h^{3/2}$.}
\label{fig:conv}
\end{figure}

The Dirichlet results agree with the rule \reff{eq:rule}.
For $m=0$ and $p\ge2$, the predicted energy rate is
$\min(p,2m+\frac32)=\frac32$, compared with measured rates of $1.54$ for
$p=2$ and $1.33$ for $p=3$.
Adding one Taylor term raises these rates to $2.22$ and $3.41$.

At fixed finite $\beta$, the improved estimate in
\Thm{\ref{thm:err}} predicts an energy rate of $\min(p,2m+2)$.
For $p=2$ and $m=0$, the predicted order $2$ agrees with the measured rates
of $2.15$ at $\beta=0$ and $2.20$ at $\beta=1$.
The optimal energy convergence in these cases therefore follows from the
analysis.  For $p=3$ and $m=0$, the predicted rate is again $2$, while the
measured rates are $1.81$ and $1.85$ on the tested meshes.
The distinction between finite $\beta$ and the Dirichlet limit comes from
the complementary weights in \Rem{\ref{rem:where}}, which confine the
half-order loss in the uniform estimate to the $\beta$-weighted part of
the residual.

The improved energy estimate does not imply optimal $L^2$ convergence.
For $p=2$ and $m=0$, the measured $L^2$ rates range from $1.77$ to $1.87$,
increasing with $m=1$ to between $2.96$ and $3.29$.
This agrees with the need for a correction discussed in
\Rem{\ref{pro:l2}}.
At $\beta=0$, where only the flux transport $\Lm$ is active, increasing
$m$ from $0$ to $1$ raises the $L^2$ rate from $1.77$ to $2.96$.

\begin{figure}[!htbp]\centering
\includegraphics[width=0.42\textwidth]{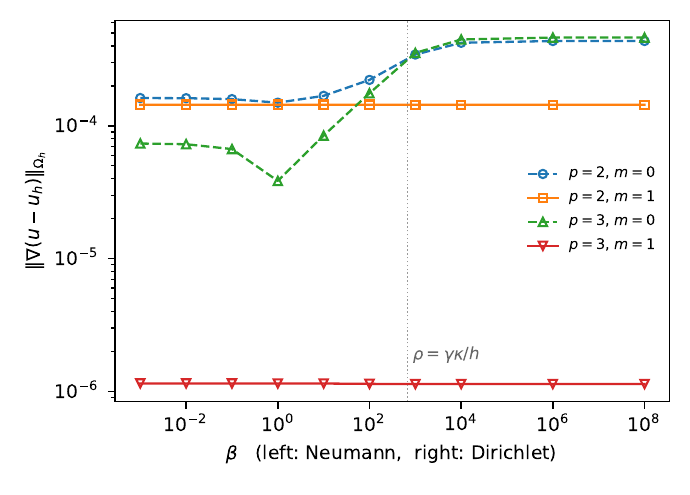}\hfill
\includegraphics[width=0.42\textwidth]{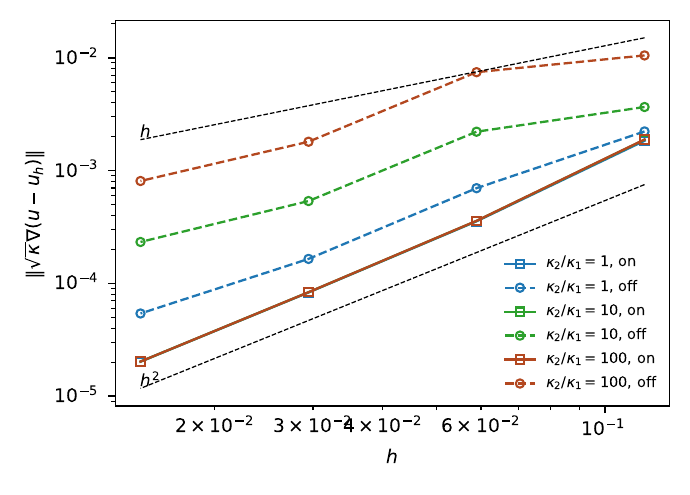}
\caption{Left: the $H^1$ error at $N=48$ as $\beta$ varies from the Neumann to the Dirichlet limit.  The curves are nearly flat for $m=1$.  Right: the interface problem with (solid) and without (dashed) the flux-jump term of \reff{eq:iAh}, for three contrasts.}
\label{fig:sweep}
\end{figure}

We also vary $\beta$ over ten orders of magnitude at $N=48$, with the
results shown in \Fig{\ref{fig:sweep}}, left.
For $m=0$, the error grows monotonically by a factor $2.69$ between the
Neumann and Dirichlet limits.
The transition occurs near $\beta\simeq\rho=\gamma\kappa/h$, which is of
order $10^3$ on this mesh.
For $m=1$, the variation is at most $0.1$ per cent over the tested range,
so the correction removes most of the observed dependence on $\beta$
in this example.

On the non-convex flower, where $\rh$ changes sign along $\Gh$, the energy
rates for $p=2$ and $m=1$ are $2.24$, $2.24$ and $2.23$ for
$\beta=0,1,\infty$.
For $m=0$, the Dirichlet rate falls to $1.22$.
The corrected method therefore retains its convergence rate despite the
changing sign of $\rh$.

\subsection{The interface problem and the flux-jump correction.}
\label{ss:imis}

We use the same circular interface with $\kappa_1=1$ and $\kappa_2=10$,
and choose the exact solutions
$u_1=\sin(1.7x)\cos(1.7y)$ and
$u_2=0.3+x^2-\frac12y^2+0.4xy$.
Both the solution jump and the flux jump are non-zero, as is the
tangential flux jump in \Rem{\ref{rem:keep}}.

For perfect bonding, the $H^1$ rate is $1.00$ for $p=1$.
For $p=2$, it increases from $1.40$ with $m=0$ to $2.17$ with $m=1$, while
for $p=3$ the corresponding rates are $1.36$ and $3.32$.
These results agree with the rule \reff{eq:rule}, and the case
$\beta_\Gamma=10$ gives similar results on these meshes.
At $N=48$ and $p=2$, we also vary $\beta_\Gamma$ over nine orders of
magnitude.  With $m=1$, the error varies by $0.3$ per cent, indicating weak
dependence on the bonding stiffness in the tested range.

To isolate the flux-jump correction, we run the method with and without
the second line of \reff{eq:iAh}.
In the latter case, the flux-jump datum is imposed directly on $\Gh$.
Both runs use $p=2$, $m=1$ and perfect bonding, with the results given in
\Tab{\ref{tab:imis}}.

\begin{table}[!htbp]\centering
\begin{tabular}{lrrrr}\toprule
contrast $\kappa_2/\kappa_1$ & \multicolumn{2}{c}{with the term} & \multicolumn{2}{c}{without} \\
 & $H^1$ rate & $L^2$ rate & $H^1$ rate & $L^2$ rate \\ \midrule
1 & 2.16 & 3.24 & 1.82 & 1.88 \\
10 & 2.17 & 3.28 & 1.40 & 1.88 \\
100 & 2.17 & 3.27 & 1.31 & 1.88 \\
\bottomrule
\end{tabular}
\caption{With and without the flux-jump term of \reff{eq:iAh}, $p=2$, $m=1$, perfect bonding.}
\label{tab:imis}
\end{table}

Omitting the flux-jump correction reduces the rate at each contrast.
At contrast $1$, the $H^1$ rate falls from $2.16$ to $1.82$, while at
contrasts $10$ and $100$ it falls from $2.17$ to $1.40$ and $1.31$,
respectively.
The tangential flux jump is
$(\kappa_1-\kappa_2)\p_{\V t}u_1+\kappa_2\p_{\V t}\jm{u}$, which need not
vanish at unit contrast.
In this example, it remains non-zero because the jump varies along the
interface.  Equal coefficients and a constant jump would make the term
vanish.
The corrected rates are almost unchanged over the three contrasts,
supporting the weighted estimate in \Thm{\ref{thm:ierr}} for these cases.
This correction adds two terms on the cut elements.

\subsection{The symmetric preconditioner.}
\label{ss:prec}

\Cor{\ref{cor:precond}} gives a contraction factor of order $h^{1/2}$ when
$A_h^0$ preconditions the corrected form, a prediction we test with $m=1$
in \Tab{\ref{tab:prec}}.
The observed factor is obtained from the iteration
$e\mapsto e-(A_h^0)^{-1}A_he$ on the homogeneous problem, starting from a
random vector.
We also report iteration counts for Richardson iteration and GMRES, both
using the same factorisation of $A_h^0$ and a stopping tolerance of
$10^{-8}$ for the relative residual.

\begin{table}[!htbp]\centering
\begin{tabular}{lrrrrr}\toprule
 & \multicolumn{2}{c}{$q$ observed} & \multicolumn{2}{c}{iterations (max)} & \\
case & $N=12$ & $N=96$ & Richardson & GMRES & $\|u_{\rm it}-u_{\rm dir}\|/\|u_{\rm dir}\|$ \\ \midrule
$\beta=0$, $p=2$ & 0.0206 & 0.0041 & 5 & 4 & $1.26\cdot10^{-9}$ \\
$\beta=0$, $p=3$ & 0.0226 & 0.0069 & 5 & 5 & $3.91\cdot10^{-10}$ \\
$\beta=1$, $p=2$ & 0.0181 & 0.0037 & 5 & 4 & $2.44\cdot10^{-10}$ \\
$\beta=1$, $p=3$ & 0.0222 & 0.0073 & 5 & 4 & $7.44\cdot10^{-11}$ \\
$\beta=\infty$, $p=2$ & 0.0186 & 0.0041 & 4 & 4 & $3.13\cdot10^{-9}$ \\
$\beta=\infty$, $p=3$ & 0.0238 & 0.0044 & 4 & 4 & $3.32\cdot10^{-9}$ \\
\bottomrule
\end{tabular}
\caption{The uncorrected symmetric form as a preconditioner, $m=1$: observed contraction factor on the coarsest and finest meshes, the largest iteration count over $N=12,\dots,96$ to a relative residual $10^{-8}$, and the largest discrepancy with the direct solution of the non-symmetric system.}
\label{tab:prec}
\end{table}

The contraction factor never exceeds $0.028$, and its smallest value on
the finest mesh is $0.0037$.
Over the four meshes, the least-squares rates range from $0.58$ to $0.88$
across the six cases, consistent with the $O(h^{1/2})$ bound of
\Cor{\ref{cor:precond}}.
At the Dirichlet limit, $\rho s\to0$, so the term responsible for the
half-order bound in \Pro{\ref{pro:pert}} vanishes from \reff{eq:Bh},
leaving an $O(h)$ perturbation.

No more than five iterations are needed in these tests, and the counts
decrease under refinement.
The last column compares the iterative and direct solutions of the
non-symmetric system, with small discrepancies consistent with the
stopping tolerance and finite precision arithmetic.
The method therefore uses one symmetric factorisation and applies the
geometric correction at each iteration.


\Section{Conclusions.}
\label{sec:concl}

We have presented a cut finite element method with boundary value
correction for Robin and Neumann conditions, including an extension to
interfaces with imperfect bonding.  By transporting the solution and its
gradient to the exact boundary and evaluating the flux with the exact
normal, the method retains a piecewise affine geometry for all integration.

The energy analysis gives a geometric error of
$O(h^{-1/2}\delh^{m+1})$ uniformly through the Dirichlet limit, but the
continuous trace inequality sharpens this to $O(\delh^{m+1})$ for fixed
finite $\beta$.  The distinction comes from the trace residual, which
carries the Dirichlet weight, and leads to the sufficient conditions
$2m+\frac32\ge p$ and $2m+2\ge p$ when $\delh=O(h^2)$.
In particular, $p=2$ and $m=0$ give second-order energy convergence for
fixed finite $\beta$, explaining the optimal energy rates in
\Tab{\ref{tab:conv}} under the stated assumptions.  The analysis does
not establish the larger slopes measured on these meshes, and the
corresponding $L^2$ estimate still requires the adjoint-consistency
argument discussed in \Rem{\ref{pro:l2}}.

For interface problems, the additional flux-jump correction accounts for
the different normals used by integration by parts and by the prescribed
jump condition.  The numerical examples show a substantial loss of
accuracy when this term is omitted, whereas the corrected method gives
energy rates close to two for quadratic elements at all three tested
contrasts.  Uniform stability remains subject to the Poincar\'e bounds
stated in the analysis, with a normalization for any uncoupled constant
mode at zero bonding stiffness.

The perturbation estimate also suggests a simple preconditioner based on
the symmetric form assembled on the affine geometry.  The preconditioned solvers require at most five iterations in the
numerical examples, with the geometric correction applied without a
separate factorisation.

\paragraph{Funding.}
EB was supported in part by EPSRC grants EP/X042650/1 and EP/V050400/1. PH
was supported by the Swedish Research Council under grant 2022-03908. MGL was supported in
part by the Swedish Research Council under grants 2021-04925 and 2025-05562, the Knut and Alice
Wallenberg Foundation under grant KAW 2025.0277, and the Swedish Research Programme Essence.

%

\end{document}